\documentclass[12pt]{article}

\usepackage{fullpage,amsmath,amsfonts}
\usepackage[all]{xy}
\usepackage{amssymb,amsthm}
\usepackage{hyperref}

\newcommand{\Set}{\mathbf{Set}}

\newcommand{\calB}{\ensuremath{\mathcal{B}}}
\newcommand{\calC}{\ensuremath{\mathcal{C}}}
\newcommand{\calD}{\ensuremath{\mathcal{D}}}
\newcommand{\calE}{\ensuremath{\mathcal{E}}}
\newcommand{\calF}{\ensuremath{\mathcal{F}}}
\newcommand{\calG}{\ensuremath{\mathcal{G}}}
\newcommand{\calI}{\ensuremath{\mathcal{I}}}
\newcommand{\calJ}{\mathcal{J}}
\newcommand{\calL}{\ensuremath{\mathcal{L}}}
\newcommand{\calS}{\ensuremath{\mathcal{S}}}

\newcommand{\opCat}[1]{\ensuremath{#1^\mathrm{op}}}

\newcommand{\Nat}{\ensuremath{\mathbb{N}}}

\newcommand{\Psh}[1]{\widehat{#1}}

\newcommand{\rmD}{\mathrm{D}}   
\newcommand{\rmE}{\mathrm{E}}

\newcommand{\aE}{\acute{\rmE}}

\newcommand{\IBD}[1]{\mathrm{IBD}_{#1}}

\newcommand{\height}{\mathfrak{h}}

\newcommand{\Sk}{\mathrm{Sk}}

\newcommand{\emin}[1]{{#1}^{\ast}}

\DeclareMathOperator{\petit}{\mathrm{P}}

\theoremstyle{plain}
\newtheorem{theorem}{Theorem}[section]
\newtheorem{corollary}[theorem]{Corollary}
\newtheorem{proposition}[theorem]{Proposition}
\newtheorem{lemma}[theorem]{Lemma}

\theoremstyle{definition}
\newtheorem{definition}[theorem]{Definition}
\newtheorem{example}[theorem]{Example}
\newtheorem{remark}[theorem]{Remark}

\title{The étendue of a combinatorial space \\ and its dimension}
\author{M. Menni\footnote{Conicet y Centro de Matemática de La Plata, Argentina}\footnote{This project has received funding from the European Union’s Horizon 2020 research and innovation programme under the Marie Skłodowska-Curie grant agreement No. 101007627, and also from CONICET (Argentina), PIP 11220200100912CO.}}
\date{\today}

\begin{document}

\maketitle

\begin{abstract}
To each simplicial set $X$ we naturally assign an étendue ${\aE X}$ whose internal logic  captures   information about the geometry of  $X$. 
In particular, we show that, for   `non-singular' objects $X$ and $Y$, the étendues  ${\aE X}$ and ${\aE Y}$ are equivalent if, and only if,   $X$ and $Y$ have the same dimension. Many of the results apply to   presheaf toposes over `well-founded' sites.
\end{abstract}

\tableofcontents

\section{Introduction}
\label{SecIntro}

Motivated by one of Lawvere's suggestions in \cite[Section II]{Lawvere91} we show that, for a certain type of gros topos $\calE$, the dimension of a `non singular'  object $X$ in $\calE$ is determined by the petit topos of $X$. The main result applies to the topos of simplicial sets \cite{GabrielZisman}, the classifier of non-trivial Boolean algebras \cite{Lawvere2021}, the topos of Ball complexes \cite{Roy2021} and  other toposes of combinatorial spaces. 
To a certain extent, the result also applies to  some presheaf toposes that are not quite combinatorial, in the sense that they do not have  locally finite sites. For instance,  the classifier of complex algebras with exactly two idempotents (i.e. the Gaeta topos determined by the field $\mathbb{C}$) \cite{MMlevelEpsilon}, which may be viewed as a first approximation to the gros Zariski topos (over the same field)  sometimes considered in algebraic geometry. Similarly, some of the results also apply to the Gaeta topos determined by the initial distributive lattice instead of the field $\mathbb{C}$ and to the Gaeta toposes determined by the theories of distributive lattices,  of  integral rigs \cite{Menni2021a} and of  positive rigs \cite{Menni+2024+Postive+Rigs}.  We will nevetherless concentrate on combinatorial examples.

In this introduction we recall some of the material in \cite[Sections~I and II]{Lawvere91} in order to  motivate and state our main result more precisely.

Although not yet completely understood, there is a recognized distinction  between  `petit' and `gros'  toposes; see  \cite[Exercice {4.10.6}]{SGA4}, \cite{CosteMichon1981}, \cite{Lawvere86}, \cite{Street2000}, \cite{Johnstone2012}, \cite{Menni2022} and also  \cite{Lawvere89} where the properties defining  \'etendues and QD toposes are considered as typical of  `petit' toposes. (We will recall the definition of {\em étendue}  just before Definition~\ref{DefLevelÉ}. These toposes are also called  {\em locally localic}. 
 A topos is {\em QD} if every object in it is a quotient of a decidable object. QD toposes are sometimes called {\em locally decidable}. See  \cite[C5.2.4 and C5.4.1]{elephant} and the references therein.) 

 Partially inspired by distinction between the classes of `gros' and `petit' toposes,  \cite[Section I]{Lawvere91} proposes a tentative clarification of the distinction  and explores the question of how one class arises from the other. In particular, given an object $X$ in a topos `of spaces',  what is the reasonable topos of pseudo-classical sheaves on $X$?  Lawvere explicitly says that he has not succeeded to give a site-invariant description but, in order to clarify the problem, he proposes to study a particular class of examples. More explicitly, given a small extensive category ${\calC}$, let ${\mathrm{G}\calC}$ be the {\em Gaeta} topos of $\calC$; in other words, the topos of finite-product preserving functors ${\opCat{\calC}\rightarrow \Set}$. For any object $X$ in ${\mathrm{G}\calC}$,  Lawvere constructs a QD-subtopos ${\mathrm{P}(X)}$ of the slice ${\mathrm{G}\calC/X}$ and formulates a conjecture relating ${\mathrm{P}(X)}$ and the dimension of $X$.

The foundations of Lawvere's  dimension theory are fully described in  \cite[Section~{II}]{Lawvere91}.
We quickly review here  the main definitions.
A subtopos ${l : \calL \rightarrow \calE}$ is {\em essential} if the inverse image functor $l^*$ has a left adjoint, denoted by ${l_! : \calL \rightarrow \calE}$. In the context we are discussing, essential subtoposes are usually  called {\em levels}, and the suggestion in \cite{Lawvere91} is to identify levels with dimensions.
More specifically, we think of ${l_! : \calL \rightarrow \calE}$ as the full subcategory of objects of dimension $l$. So, for an object $X$ in $\calE$, we write ${\dim X \leq l}$ if the counit ${l_! (l^* X) \rightarrow X}$ is an isomorphism; in this case we may also say that $X$ is $l$-skeletal.  It is relevant to mention that all this is fully consistent with the standard dimension theory of simplicial sets; see \cite{GabrielZisman, KennettEtAl, Menni2019a, MenniPAMS} and references therein.

We now quote from \cite[Section~{II}]{Lawvere91}:
``The idea behind the identification of the levels in a category of
Being with dimensions is that a higher level is a more determinate
general Becoming, that is, it contains spaces having in them possibly-more-varied information for determining processes. Thus one
conjectures that ${\dim X}$ only depends on the category ${\mathrm{P}(X)}$ of particular
Becoming associated to $X$ (and not on the important structure sheaf
which recalls for the little category the big environment in which it
was born). In other words, if we have an equivalence of categories
${\mathrm{P}(X) \cong \mathrm{P}(Y)}$, then $X,Y$ should belong to the same class of UIO levels
within the category of Being in which they are objects. Suitable hypotheses
to make this conjecture true should begin to clarify the relationships
between the two suggested philosophical guides."

We find hypotheses to make a variant of the above conjecture true.

Roughly speaking, we use `locally localic' instead of `locally decidable' as a choice to model particular Becoming.
This choice allows us  not only to prove the conjecture, but also to give a site-invariant description of independent interest.
In fact, we may use this description to outline the results in more detail.

Recall from \cite[Definition~{A4.6.1}]{elephant} that a geometric morphism ${f : \calF \rightarrow \calS}$ is {\em localic} if every object of $\calF$ is a subquotient of one of the form ${f^* A}$ for some object $A$ in $\calS$. 
The concept of étendue is usually formulated with a tacit base. Let us be more explicit.
We say that a geometric morphism ${g : \calG \rightarrow \calS}$ is  an {\em étendue} if there is a well-supported object $G$ in $\calG$ such that the composite 
\[\xymatrix{
 \calG/G \ar[r] & \calG \ar[r]^-g & \calS
}\]
is localic, where ${\calG/G \rightarrow \calG}$ is the canonical geometric morphism induced by ${G}$.

\begin{definition}[Level~{é}]\label{DefLevelÉ}
A geometric morphism ${f : \calF \rightarrow \calS}$ has a {\em level~{é}} if $\calF$ has a largest level ${\calL \rightarrow \calF}$ such that ${\calL \rightarrow \calF \rightarrow \calS}$ is an étendue.
\end{definition}

The concept of Level~{é} seems to be new. 
The definition is analogous to that of level~{$\epsilon$} discussed in \cite{MMlevelEpsilon} which is also defined as the largest level with a certain property.
Anyway, we will see that many  geometric morphisms  have a non-trivial level~{é} and this will allow us to model particular Becoming.

To motivate the idea we quote again from  \cite{Lawvere91}: ``A general category of Being, particular categories of Becoming: this is a suggested philosophical guide for sorting the two original kinds of toposes and what they have become. The unity and cohesiveness of Being provides the basis for Becoming, and the historicity and controlled variability of Becoming produces new Being from old."   

For this reason we are interested in  geometric morphisms ${p : \calE \rightarrow \calS}$ such that, for every object $X$ in $\calE$, the composite ${\calE/X \rightarrow \calE\rightarrow \calS}$  has a level~é, which will be denoted by ${\aE X \rightarrow \calE/X}$.
For such a geometric morphism $p$, the topos $\calE$ is to be thought of as the general category of Being, and the étendues ${\aE X}$ as the particular categories of Becoming, one for each object $X$ in $\calE$.

At the risk of being somewhat vague, we make an informal comment as an attempt to give some additional intuition about the toposes ${\aE X}$.
Even in well-understood toposes of spaces,  it is not clear what is an `open' subobject of  a space $X$. 
For instance, what is an open subobject of a simplicial set? The vague suggestion  is that the topos ${\aE X}$ is a surrogate for the {\em a priori} non-existent  locale of open subobjects of $X$. The topos ${\aE X}$ is only locally localic, but we will see that, in many cases, it does capture the dimension of $X$.  In particular,   we will show that,  for every pair of `non-singular' simplicial sets $X$, $Y$,   ${\aE X \cong \aE Y}$ if, and only if,  $X$ and $Y$  have the same dimension.

The way in which  ${\aE X}$ determines the dimension of $X$ is inspired by  \cite{MarraEtAl2018} where, among other things, the authors show  that, for a compact polyhedron $P$, the Heyting algebra of open subpolyhedra of $P$ captures the  dimension of $P$ through the bounded-depth formulae. The proof relies on fairly combinatorial techniques involving triangulations but it does not seem to lead naturally to a result about purely combinatorial objects such as simplicial sets. In a sense, the present paper is an attempt to arrive at such a result using \cite{Lawvere91} as a guide.

In order to give more detail, and also for brevity, we introduce the following notation for a Heyting formula that will play  an important role throughout  the paper. 
Indeed, let 
\[ \gamma(x, y) := x \vee (x \Rightarrow y) \]
in variables $x, y$.

For instance, if $H$ is a biHeyting algebra  and ${a, b \in H}$ then, ${\top \leq \gamma(a,b)}$ if and only if ${\partial a \leq b}$, that is, if and only if, the coHeyting boundary of $a$  is below $b$.

Also, using the internal language of a topos with subobject classifier $\Omega$,  let 
\[
\begin{array}{lclcl}
 \IBD{-\infty} & := & \bot \\
 \IBD{0} & := &  (\forall x : \Omega)(\gamma(x, \IBD{-\infty})) & = &  (\forall x : \Omega)(x \vee \neg x) \\
     \IBD{n+1} & := & (\forall x:\Omega)(\gamma(x, \IBD{n})) \\
     \IBD{\infty} & := & \top
\end{array}
\]
for every ${n \in \Nat}$. Needless to say, these are `closed' variants of the Bounded Depth formulas used in \cite{MarraEtAl2018}. Notice that the interpretation of $\IBD{n}$ in any topos is a subterminal for every ${n \in \{-\infty, \infty \} + \Nat}$. Trivially, a topos validates ${\IBD{-\infty} = \bot}$ if and only if it is degenerate. Also, it is well-known that a topos validates ${\IBD{0}}$  if and only if it is Boolean \cite[Section~{VI.5}]{maclane2}.

As a corollary of our results we will show that, for every `non-singular'  simplicial set $X$, and for every ${n \in \{-\infty, \infty\} + \Nat}$, ${\dim X \leq n}$ (i.e. $X$ is $n$-skeletal) if, and only if, ${\IBD{n}}$ holds in ${\aE X}$. In other words, the dimension of $X$ depends only on  ${\aE X}$ (and not on the composite ${\aE X \rightarrow \Psh{\Delta}/X \rightarrow \Psh{\Delta}}$ which recalls for ${\aE X}$ ``the big environment in which it was born").

In Section~\ref{SecWidespread} we recall the notion of widespread subterminal \cite{Johnstone1979au} and identify when a map with terminal codomain factors through  ${(\forall x : \Omega)(\gamma(x, u))}$ where $u$ is a point of $\Omega$. This is applied in  Section~\ref{SecBoundaries} to characterize when ${\IBD{n}}$ holds in presheaf toposes. In Section~\ref{SecNonSingular} we introduce minimal objects. These provide a generalization of the notion   of  non-degenerate simplex  \cite[Section~{II.3}]{GabrielZisman} and are  used to define a notion of strongly regular presheaf. Such presheaves are,  intuitively,  `sufficiently non-singular'. To give a more precise  idea here we simply mention that the category of strongly regular simplicial sets coincides with that studied  in \cite{Zisman2009} and it contains the category of non-singular simplicial sets considered in  \cite{WaldhausenEtAl}.
Minimal objects are also used in Section~\ref{SecD} to assign, to each presheaf $X$, an étendue  ${\rmD X}$. 
We give some examples and, in particular, observe that, for non-singular $X$, ${\rmD X}$ is localic.

In Section~\ref{SecLevelE} we prove that, for certain small $\calC$, ${\rmD X}$ coincides with ${\aE X}$ for every $X$ in $\Psh{\calC}$.
As a corollary,  ${\aE X}$ exists for every simplicial set $X$.

In Section~\ref{SecBeingAndBecoming} we show that,  if  ${\dim X \leq n}$ then $\IBD{n}$ is true in ${\rmD X}$ and, in Section~\ref{SecWF} we prove that for `well founded' sites,  the converse holds for strongly regular $X$.

We end the present section with a couple of comments about notation.

For a presheaf $X$ on a small  $\calC$ we denote its category of elements by ${\calC/X}$.
This is the notation used in \cite{Lawvere91} and it is convenient, in particular, because the explicit description of slices of presheaf toposes (see e.g. \cite[Proposition~{A1.1.7}]{elephant}) admits an efficient and informative formulation as follows
\[ \Psh{\calC}/X \cong \Psh{\calC/X} \]
for any $X$ in $\Psh{\calC}$. The objects of ${\calC/X}$ will be denoted by ${(x, C)}$ with ${x \in X C}$.

As usual,   $\Delta$ denotes the category of ordered sets ${[n] = \{ 0 < \ldots < n \}}$ with ${n \in \Nat}$ and monotone maps between them, so that ${\Psh{\Delta}}$ is the topos of simplicial sets. This will be one of the main examples in the paper.
The levels of ${\Psh{\Delta}}$ are those induced by the truncations ${\Delta_n \rightarrow \Delta}$ where, as is customary,  we let ${\Delta_n \rightarrow \Delta}$ be the full subcategory determined by the objects ${[m]}$ with ${m \leq n \in \Nat}$.
(See \cite[Corollary~{3.3}]{MenniPAMS} and also \cite[Proposition~{2.5}]{Menni2019a}.)

Also, let $\mathbb{F}$ be the category of finite non-empty sets and functions.
It shares many properties with $\Delta$ but it is not a Reedy category.
The topos ${\Psh{\mathbb{F}}}$ and its possible applications in Combinatorial Topology are discussed in \cite[Section~1]{Lawvere2021}.

\section{Widespread subterminals}
\label{SecWidespread}

Following \cite[D4.5.9]{elephant} we say that an element $w$  in a distributive lattice $H$ is {\em widespread} if the lattice ${\{ v \in H \mid w \leq v\}}$ is complemented.
For instance, the bottom element is widespread if and only if the lattice is a Boolean algebra.

Now recall our notation 
\[ \gamma(x, y) := x \vee (x \Rightarrow y) \]
in variables $x, y$, and consider the following characterization of widespread elements.

\begin{lemma}\label{LemCharWidespread}
If $H$ is a Heyting algebra then, for any ${w \in H}$, the following are equivalent:
\begin{enumerate}
    \item $w$ is widespread.
    \item For every ${v \in H}$, ${\top \leq \gamma(v,w) = v \vee (v\Rightarrow w)}$.
\end{enumerate}
\end{lemma}
\begin{proof}
Assume first that $w$ is widespread and let ${v \in H}$.
Then ${v \vee w}$ has a complement $v'$ in the lattice of elements above $w$.
That is, ${\top \leq v \vee v'}$ and ${v \wedge v' = w}$.
Hence, 
\[ \top \leq  v \vee v' \leq  v \vee (v \Rightarrow w) \]
as we needed to show. The converse is left for the reader.
\end{proof}

For instance, if $\calE$ is a topos and  $l$ is a level of $\calE$ with monic skeleta then,  an object $X$ in $\calE$ {\em has $l$-skeletal boundaries} in the sense of \cite[Definition~{3.2}]{Menni2019a} if the $l$-skeleton of $X$ is widespread as a subobject of $X$.

Also in a topos, the interpretation of the formula below in context ${y : \Omega}$
\[ (\forall x : \Omega)( \gamma(x, y) )  \]
is a subobject of $\Omega$ and its domain is  called the {\em Higgs object} of the  topos \cite[D4.5.9]{elephant}.

If ${u : 1 \rightarrow \Omega}$ is a point of $\Omega$ in a topos, then the interpretation of ${ \gamma(x, u) }$
in context ${x : \Omega}$ is a subobject of $\Omega$ and so, the interpretation of 
\[ (\forall x: \Omega)( \gamma(x, u) )
\]
is a subterminal.

\begin{lemma}\label{LemForDefInternallyWidespread}
Let $U$ be subterminal in a topos and let ${u : 1 \rightarrow \Omega}$ be the corresponding characteristic map. Then, for any object $X$, the unique map ${X \rightarrow 1}$ factors through
the subterminal 
\[ (\forall x: \Omega)( \gamma(x, u) )
\]
if, and only if, for every map ${Y \rightarrow X}$, the subobject  ${Y \times U \rightarrow Y}$ is widespread in the Heyting algebra of subobjects of $Y$.
\end{lemma}
\begin{proof}
This is a simple exercise in the internal language.
There is a map from $X$ to the indicated subterminal iff, for every ${e : Y \rightarrow X}$ and   every ${\chi : Y \rightarrow \Omega}$,
the formula ${\gamma(\chi, u)}$ holds. Equivalently, for every subobject  $V$ of $Y$,
\[ \top \leq V \vee (V \Rightarrow (Y \times U)) \]
in the Heyting algebra of subobjects of $Y$.
 So the equivalence follows from Lemma~\ref{LemCharWidespread}.
\end{proof}

As an extreme case we obtain the following.

\begin{lemma}\label{LemDefInternallyWidespread}
Let $U$ be subterminal in a topos and let ${u : 1 \rightarrow \Omega}$ be the corresponding characteristic map. Then the following are equivalent.
\begin{enumerate}
    \item For any $X$, ${X \times U \rightarrow X}$ is widespread in the Heyting algebra of subobjects of $X$.
    \item The sentence ${  (\forall x : \Omega)(\gamma(x, u) )}$ holds in the internal logic of the topos.
    \item The closed subtopos determined by $U$ is Boolean.
\end{enumerate}
\end{lemma}
\begin{proof}
The equivalence between the first two items follows from Lemma~\ref{LemForDefInternallyWidespread}.
The equivalence with the third item follows from \cite[Lemma~4]{Johnstone1979au}.
\end{proof}

If the equivalent conditions of Lemma~\ref{LemDefInternallyWidespread} hold then  $U$ is said to be {\em (internally) widespread} \cite[Section~{D4.5}]{elephant}.

\section{Boundaries in presheaf toposes}
\label{SecBoundaries}

In this section we characterize when $\IBD{n}$ holds in a presheaf topos.
Fix a small category $\calD$ and the corresponding presheaf topos ${\Psh{\calD}}$.
If we fix a subobject ${w : W \rightarrow X}$ in ${\Psh{\calD}}$ then, for any subobject ${u : U \rightarrow X}$, the subobject ${u \Rightarrow w : (U \Rightarrow W) \rightarrow X}$ may be described as follows:
\[ (U \Rightarrow W) D = \{ x \in X D \mid \textnormal{ for all } f : E \rightarrow D, x\cdot f \in U E \textnormal{ implies } x \cdot f \in W E\} \subseteq X D \]
for each object $D$ in $\calD$.

The next result is a mild generalization of \cite[Lemma~{5.1}]{MenniPAMS}.
For each $K$ in $\calD$ and ${k \in X K}$, the image of the corresponding ${k : \calD(-, K) \rightarrow X}$ will be denoted by ${\xi_k : \Xi_k \rightarrow X}$ so that, for every $J$ in $\calD$, ${j \in \Xi_k J \subseteq X J}$ if and only if there is a ${g :J \rightarrow K}$ such that ${k \cdot g = j}$.

\begin{lemma}\label{Lem5.1again}
Let  ${w : W \rightarrow X}$ be a subobject in ${\Psh{\calD}}$.
For any $D$ in $\calD$ and any ${x : \calD(-,D) \rightarrow X}$, the following are equivalent:
\begin{enumerate}
\item For any subobject ${v : V \rightarrow X}$, $x$ factors through the subobject ${\gamma(v, w)}$  of $X$.
\item For every ${f : E \rightarrow D}$  in $\calD$, $x$ factors through  the subobject ${\gamma(\xi_{x \cdot f}, w)}$ of $X$.
\item For every ${f : E \rightarrow D}$ in $\calD$, either ${x\cdot f \in W E}$ or there is a ${g : D \rightarrow E}$ such that ${x \cdot (f g) = x}$.
\end{enumerate}
\end{lemma}
\begin{proof}
The first item easily implies the second.
To prove that the second implies the third let ${f : E \rightarrow D}$ and consider the subobject ${\xi_{x \cdot f} : \Xi_{x \cdot f} \rightarrow X}$. By hypothesis, ${x \in \Xi_{x\cdot f} D}$  or ${x \in (\Xi_{x \cdot f} \Rightarrow W) D}$. If ${x \in \Xi_{x\cdot f} D}$ then there is a ${g : D \rightarrow E}$ such that ${(x \cdot f) \cdot g = x}$. If ${x \in (\Xi_{x \cdot f} \Rightarrow W) D}$ then, for all ${h : B \rightarrow D}$, ${x \cdot h \in \Xi_{x \cdot f} B}$ implies ${x \cdot h \in W B}$. Taking ${h = f : E \rightarrow D}$ we conclude that ${x\cdot f \in W E}$.

To prove that the third item implies the first, let ${v : V\rightarrow X}$ be a subobject of $X$ and assume that ${x \not\in V D \subseteq X D}$. We show that ${x \in (V \Rightarrow W) D \subseteq X D}$. To do this let ${f : E \rightarrow D}$ be such that ${x \cdot f \in V E}$. Then, by hypothesis, either ${x \cdot f \in W E}$ or there is a ${g : D \rightarrow E}$ such that ${x \cdot (f g) = x}$. If such a $g$ exists then, as ${x\cdot f \in V E}$, ${(x\cdot f) \cdot g \in V D}$. So ${x \in V D}$, which contradicts our assumption.
Hence, ${x \cdot f \in V E}$ implies ${x\cdot f \in W E}$. That is ${x \in (V \Rightarrow W) D}$ as we needed to show.
\end{proof}

The following is the obvious variant of \cite[Lemma~{5.2}]{MenniPAMS}.

\begin{lemma}\label{Lem5.2again}
For any map ${e : E \rightarrow D}$ in $\calD$  and any subobject ${w : W \rightarrow \calD(-,D)}$ in $\Psh{\calD}$ the following hold:
\begin{enumerate}
\item   The map ${\calD(-, e) : \calD(-, E) \rightarrow \calD(-,D)}$  factors through the subobject ${v \vee (v \Rightarrow w)}$ for any subobject ${v : V \rightarrow \calD(-,D)}$ in $\Psh{\calD}$, if and only if, for every ${f : F \rightarrow E}$ in $\calD$, either ${e f \in W F}$ or there is a ${g : E \rightarrow F}$ such that ${e f g = e}$.
\item If $e$ is monic, the above is equivalent to: for every ${f : F \rightarrow E}$ in $\calD$, either ${e f \in W F}$ or $f$ has a section.
\end{enumerate}
\end{lemma}
\begin{proof}
The first item follows from  Lemma~\ref{Lem5.1again}, and the second item follows from the first.
\end{proof}

Taking $e$ to be the identity on $D$ we deduce the characterization of widespread sieves below.

\begin{lemma}\label{LemWidespreadSieve}
A subobject ${w : W \rightarrow \calD(-,D)}$ in $\Psh{\calD}$ is widespread if, and only if, for every ${f : F \rightarrow D}$, ${f \in W F}$ or $f$ has a section.    
\end{lemma}

Subterminals in $\Psh{\calD}$ correspond to sieves on $\calD$.
That is, subsets of objects in $\calD$ such that for every ${E \rightarrow D}$ in $\calD$, if $D$ is in the subset then so is $E$.

\begin{proposition}\label{PropMeaningOfIBB}
Let ${U}$ be a subterminal in ${\Psh{\calD}}$ and let ${u : 1 \rightarrow \Omega}$ be its classifying map. Then the interpretation of 
\[   (\forall v : \Omega)(\gamma(v,u) ) \]
corresponds to the sieve consisting of the objects $D$ such that,  for every ${c : C \rightarrow D}$, $C$ is in $U$ or  $c$ is an iso.
\end{proposition}
\begin{proof}
Let ${S \rightarrow 1}$ be the interpretation of the formula in the statement.
For any $D$ in $\calD$, Lemma~\ref{LemForDefInternallyWidespread} implies that:  
the unique map ${!_D : \calD(-,D) \rightarrow 1}$ factors through  ${S \rightarrow 1}$  
if and only if for every ${e : E \rightarrow D}$, 
${W_E = \calD(-,E) \times U \rightarrow \calD(-,E)}$ is widespread in the poset of subobjects of ${\calD(-,E)}$. Notice that $W_E$ is the sieve of maps with codomain $E$ whose domain is in $U$.
So, by Lemma~\ref{LemWidespreadSieve}, we get that:
\begin{itemize}
\item $!_D$ factors through ${S \rightarrow 1}$ if, and only if, for every ${e : E \rightarrow D}$ and ${f : F \rightarrow E}$, $F$ is in $U$  or $f$ has a section.
\end{itemize}
We need to check that the highlighted condition above is equivalent to the one in the statement.
Assume first that the condition in the statement holds and let ${e : E \rightarrow D}$ and ${f : F \rightarrow E}$. By hypothesis, taking ${c = e}$ we obtain that $E$ is in $U$ or $e$ is an iso. 
Taking ${c = e f}$ we deduce that $F$ is in $U$ or $e f$ is an iso.
So, if $F$ is not in $U$ then, ${e f}$ is an iso and, since $U$ is ideal, $E$ is not in $U$, so $e$ is an iso. That is, if $F$ is not in $U$, then $f$ is an iso.
Altogether, $F$ is in $U$ or $f$ is an iso.

For the converse, assume that the highlighted condition holds and let ${c : C \rightarrow D}$.
Taking ${e = id_D : D \rightarrow D}$ and ${f = c: C \rightarrow D}$ in that condition we deduce that $C$ is in $U$ or $c$ has a section, say, ${s : D \rightarrow C}$.
Taking ${e = c : C \rightarrow D}$ and ${f = s : D \rightarrow C}$ in the highlighted condition we deduce that $C$ is in $U$ or $s$ has a section. Then, $C$ is in $U$ or $c$ is an iso.
\end{proof}

In order to describe the group of automorphisms of the subobject classifier of an arbitrary presheaf topos, \cite[Example~{13}]{Johnstone1979au} introduces the following terminology: an object $E$ is called {\em extreme} if every map with domain $E$ is an isomorphism.

\begin{corollary}\label{CorPTJ-Example13}
A sieve $U$ on $\calD$ is internally widespread as a subterminal in $\Psh{\calD}$ if, and only,   every object not in $U$ is extreme.
\end{corollary}
\begin{proof}
By Proposition~\ref{PropMeaningOfIBB}, $U$ is internally widespread if, and only if, for every map ${c : C \rightarrow D}$, $C$ is in $U$ or $c$ is an iso.
\end{proof}

Corollary~\ref{CorPTJ-Example13} has a more elegant proof in  \cite[Example~{13}]{Johnstone1979au} using the relation with closed subtoposes, but we need the additional information provided by Proposition~\ref{PropMeaningOfIBB}; in particular, to  calculate the interpretation of the Bounded Depth formulas.
Recall that we defined:
\[
\begin{array}{lclcl}
 \IBD{-\infty} & := & \bot \\
 \IBD{0} & := &  (\forall x : \Omega)(\gamma(x, \IBD{-\infty})) & = & (\forall x : \Omega)(x \vee (x\Rightarrow \bot)) \\
     \IBD{n+1} & := & (\forall x:\Omega)(\gamma(x, \IBD{n})) & = & (\forall x :  \Omega)(x \vee (x\Rightarrow \IBD{n}) )
\end{array}
\]
for every ${n \in \Nat}$.
As we have already mentioned, each ${\IBD{n}}$ determines a subterminal in any topos.
In particular, in the topos of presheaves on the small category $\calD$, each ${\IBD{n}}$ determines a sieve on $\calD$.
Of course,  ${\IBD{-\infty}}$ determines the empty sieve.

\begin{corollary}\label{CorDim0}
 ${\IBD{0}}$ determines the sieve of objects $D$ such that, every map with codomain $D$ is an iso. 
\end{corollary}
\begin{proof}
By Proposition~\ref{PropMeaningOfIBB}, the sieve determined by ${\IBD{0}}$ is the sieve of objects $D$ such that for every ${f : F \rightarrow D}$, $F$ is in the empty sieve or $f$ is an isomorphism.
\end{proof}

The following is the analogue of \cite[Lemma~{2.6}]{MarraEtAl2018} that we need.

\begin{proposition}\label{PropCharIBD} For  every ${n\in\Nat}$, ${\IBD{n}}$ determines the sieve on $\calD$ consisting of the objects $D$ such that, in every sequence
\[ \xymatrix{
\cdot \ar[r] & \cdot \ar[r] & \cdots \ar[r] & \cdot \ar[r] & D
}\]
of ${n + 1}$ maps in $\calD$, one of the maps is an isomorphism. 
\end{proposition}
\begin{proof}
The case ${n = 0}$ is Corollary~\ref{CorDim0}.
Assume that the result holds for $n$ and consider the case ${n+1}$.
By Proposition~\ref{PropMeaningOfIBB}, $D$ is in the sieve determined by $\IBD{n+1}$  iff for every ${c : C \rightarrow D}$, $C$ is in the sieve determined by ${\IBD{n}}$ or $c$ is an iso.
By the inductive hypothesis, $c$ is an iso or, in every sequence as on the left below
\[ \xymatrix{
\cdot \ar[r] & \cdot \ar[r] & \cdots \ar[r] & \cdot \ar[r] & C &&&  \cdot \ar[r] & \cdot \ar[r] & \cdots \ar[r] & \cdot \ar[r] & D
}\]
with ${n + 1}$ maps in $\calD$, one of the maps is an isomorphism.
Equivalently, for every sequence   of ${n + 2}$ maps in $\calD$ as on the right above, one of the maps is an isomorphism.
\end{proof}

\section{Minimal objects and non-singular presheaves}
\label{SecNonSingular}

We recall the definitions of minimal and preterminal object, and use them to define two related notions of  `non-singular' presheaf.

\begin{definition}\label{DefMinimalObject}
An object $C$  is called {\em minimal} if every map with domain $C$ is monic.    
\end{definition}

It is relevant to stress that, in general, minimal objects are not closed under subobjects.

\begin{definition}\label{DefPreterminal}
An object $C$  is  {\em preterminal} if, for every object $B$, there is at most one map ${B \rightarrow C}$.
\end{definition}

In contrast with minimal objects we have:

\begin{remark}\label{RemSubpreterminal}
If ${m : A \rightarrow B}$ is monic and $B$ is preterminal then so is $A$.
\end{remark}

It is immediate that preterminal implies minimal; the referee made the following more subtle observation.

\begin{lemma}\label{LemReferee} 
An object $C$ is preterminal if, and only if, $C$ is minimal and, for every pair of maps ${u, v : B \rightarrow C}$ there exists an ${f : C \rightarrow D}$ such that ${f u = f v}$. 
\end{lemma}
\begin{proof}
As already mentioned, preterminal easily implies minimal but also, if $C$ is preterminal, the identity on $C$ witnesses the existence of an $f$ as in the statement. For the converse just notice that minimality of $C$ implies that $f$ is monic and therefore ${u = v}$.
\end{proof}

\begin{remark}\label{RemMinimalIffPreterminal}
It is clear from Lemma~\ref{LemReferee}  that, in the presence of terminal object, or of coequalizers, minimality and preterminality are equivalent.
In the presence of terminal object, the two notions reduce to that of subterminal object.
In particular, this is the case in slice categories.
\end{remark}

For the rest of the section, fix a small category $\calC$.

\begin{remark}[Non-degenerate figures]\label{RemNonDegenerateFigures}
For a presheaf $X$ on $\calC$, an object ${(x, C)}$ in ${\calC/X}$ is minimal if, and only if, for every ${f : C \rightarrow D}$ in $\calC$ and ${y\in X D}$ such that ${y\cdot f = x}$, $f$ is monic. 
\end{remark}

So, for instance, for a simplicial set $X$,  an object in ${\Delta/X}$  is minimal  if, and only if, it is non-degenerate in the classical  sense of \cite[Definition~{II.3.1}]{GabrielZisman}. This particular case  suggests to picture minimal objects in categories of elements as `non-degenerate' figures.
Pursuing  this intuition, the following definition is inspired by that of {\em ensemble simplicial fortement régulier} \cite[Définition~{1.3}]{Zisman2009}.

\begin{definition}\label{DefOrdinary}
An object $X$ in ${\Psh{\calC}}$ is {\em strongly regular} if, for every monic map $m$  in ${\calC/X}$  with  minimal codomain, the domain of $m$ is also minimal.
\end{definition}

Intuitively, strong regularity is a weak form of `non-singularity'.
In order to strengthen this intuition we compare strongly regular presheaves with a stronger notion also motivated by a relevant, and perhaps more familiar, type of simplicial set.

\begin{remark}\label{RemMonicFigures}
Again, for presheaf $X$ on $\calC$,  an object  ${(x, C)}$ in ${\calC/X}$ is preterminal if,  and only if,   the   corresponding  ${x: \calC(-, C) \rightarrow X}$ is monic.
\end{remark}

Hence,  it is clear that the concept introduced below is a natural generalization of the notion of {\em non-singular} simplicial set used in \cite[Definition~{1.2.2}]{WaldhausenEtAl}.

\begin{definition}\label{DefNonSingular}
 An object $X$ in $\Psh{\calC}$ is {\em non-singular} if  every minimal object of ${\calC/X}$  is preterminal.
\end{definition}

Equivalently, by Remarks~\ref{RemNonDegenerateFigures} and \ref{RemMonicFigures}, a presheaf $X$ on $\calC$ is non-singular if, for every non-degenerate figure ${(x, C)}$ of $X$, the corresponding ${x : \calC(-,C) \rightarrow X}$ is monic.

For example,  representable presheaves are non-singular because, for any object $C$ in $\calC$, ${\calC/\calC(-,C) \cong \calC/C}$ and, in slice categories, minimal objects are preterminal by Remark~\ref{RemMinimalIffPreterminal}.

\begin{proposition}\label{LemNonSingularImpliesOrdinary} Every non-singular presheaf is strongly regular.
\end{proposition}
\begin{proof}
Let $X$ be a non-singular presheaf on $\calC$. To prove that it is strongly regular, let ${m :(w,B) \rightarrow (x, C)}$ be monic in ${\calC/X}$ with minimal codomain. Then ${(x, C)}$ is preterminal because $X$ is non-singular, so ${(w, B)}$ is preterminal by Remark~\ref{RemSubpreterminal} and, hence, minimal by Lemma~\ref{LemReferee}.
\end{proof}

The category of non-singular simplicial sets  is proved to be cartesian closed in \cite{FjelboRognes2022}.

\section{Sites of minimal figures}
\label{SecD}

For a category $\calC$, the full subcategory determined by minimal objects will be denoted by  

\[ \emin{\calC} \rightarrow \calC . \]
Notice that every map in $\emin{\calC}$ is monic.

It is well-known that if $\calC$ is small then every full subcategory ${\calB \rightarrow \calC}$ induces a level ${\Psh{\calB} \rightarrow \Psh{\calC}}$ whose inverse image is restriction along the fully faithful ${\opCat{\calB} \rightarrow \opCat{\calC}}$.
See, for instance, \cite[Example~A4.5.2]{elephant} and \cite[Remark~{4.5}]{KellyLawvere89}.
In particular,  the subcategory ${\emin{\calC} \rightarrow \calC}$ of minimal objects induces such a level.

\begin{lemma}\label{LemSubEtendue}
The composite ${\Psh{\emin{\calC}} \rightarrow \Psh{\calC} \rightarrow \Set}$ is an étendue.
\end{lemma}
\begin{proof}
Since  every map in $\emin{\calC}$ is monic,  the result follows from \cite[Lemma~{C5.2.4}]{elephant}.    
\end{proof}

For any $X$ in $\calC$, the inclusion ${\emin{(\calC/X)} \rightarrow \calC/X}$ induces  level of ${\Psh{\calC}/X}$ that we denote by 
\[ \rmD X :=  \Psh{\emin{(\calC/X)}} \rightarrow \Psh{\calC/X} \cong \Psh{\calC}/X \]
and is such that  the composite ${\rmD X \rightarrow \Psh{\calC}/X \rightarrow \Psh{\calC} \rightarrow \Set}$ is an étendue by Lemma~\ref{LemSubEtendue}.

\begin{lemma}\label{LemNonSingularImpliesLocalic} 
If $X$ in $\Psh{\calC}$ is non-singular then ${\rmD X \rightarrow \Psh{\calC}/X \rightarrow \Set}$ is localic.
\end{lemma}
\begin{proof}
If $X$ is non-singular then ${\emin{(\calC/X)}}$ is easily seen to be a pre-order so the result follows from \cite[Example~{A4.6.2(d)}]{elephant}.
\end{proof}

In particular, ${\rmD X}$ is localic for every representable $X$.

\begin{example}[Strongly regular $Y$, non-localic $\rmD Y$]\label{ExNonLocalicNew} 
Let $Y$ in $\Psh{\Delta}$ be the codomain of the coequalizer ${X = \Delta(-,[1]) \rightarrow Y}$ of the two points  ${1 \rightarrow   X}$.
We picture  $Y$ as a node with a loop. The resulting site  ${\emin{(\Delta/Y)}}$  is the category
\[\xymatrix{
\bullet \ar[r]<+1ex>^-{0} \ar[r]<-1ex>_-{1} & \iota
}\]
which is obviously not a poset. In this case,  ${\rmD Y}$ is a not localic.
\end{example}

In the context of the pre-cohesive topos ${\Psh{\Delta_1}}$ of reflexive graphs, the example above is used in \cite{Lawvere86} to illustrate the distinction between toposes `of spaces' and `generalized spaces'. As a corollary of the results there, we may conclude that, for the reflexive graph $Y$ constructed as in Example~\ref{ExNonLocalicNew}, the topos ${\rmD Y}$, which may be identified with the topos  of non-reflexive graphs,  coincides with the category  of maps over $Y$ that have discrete fibers.

\begin{example}[Non strongly-regular $Z$, localic $\rmD Z$]
The three monomorphisms ${[1] \rightarrow [2]}$ in $\Delta$  induce subobjects ${\Delta(-,[1]) \rightarrow \Delta(-,[2])}$ in ${\Psh{\Delta}}$  that we can picture as the edges of a solid triangle. Denote their join by ${B\rightarrow \Delta(-,[2])}$. 
If we take the pushout
\[\xymatrix{
B \ar[d] \ar[r]^-{!} & 1 \ar[d]^-{!} \\
\Delta(-,[2]) \ar[r] & Z
}\]
then $Z$ has a non degenerate 2-simplex and a non degenerate 0-simplex, but a degenerate 1-simplex.
It follows that ${\emin{(\Delta/Z)}}$ is the total order with two elements.
Notice that, in this case, ${\rmD Z}$ satisfies $\IBD{1}$ by Proposition~\ref{PropCharIBD}, although $Z$ is 2-dimensional in the classical sense.
See also \cite[Example~{1.5}]{Menni2019a}.
\end{example}

We leave it as an exercise to study the case of the pre-cohesive topos $\calE$ of Ball complexes \cite{Roy2021} and discover if the petit topos of globular sets introduced in \cite{Street2000} is of the form ${\rmD X}$ for some object $X$ in $\calE$. Notice that, for a negative answer, it would be enough to show that the topos is not an étendue.

\section{Level~{é}}
\label{SecLevelE}

In this section we relate, for certain presheaf toposes over $\Set$,  the étendues   introduced in Section~\ref{SecD} and the  level~{é} introduced in Definition~\ref{DefLevelÉ} which, roughly speaking, is the largest level which is an étendue over the given base.

In general, we do not know if ${\Psh{\emin{\calC}} \rightarrow \Psh{\calC}}$ is level~{é} of ${\Psh{\calC} \rightarrow \Set}$, but it is related because, as the next result shows, under certain conditions, it is the largest among the   `étendue levels'  induced by full subcategories.

\begin{lemma}\label{LemPresheafSubetenduesAreBelowTheCanonicalOne}
Let $\calC$ be a small category such that every map  factors as a split-epimorphism followed by a monomorphism.
Let ${\calB \rightarrow \calC}$ a full subcategory.
If the composite ${\Psh{\calB} \rightarrow \Psh{\calC} \rightarrow \Set}$ is an étendue then 
${\calB \rightarrow \calC}$ factors through ${\emin{\calC} \rightarrow \calC}$; and so ${\Psh{\calB} \rightarrow \Psh{\calC}}$ factors through ${\Psh{\emin{\calC}} \rightarrow \Psh{\calC}}$.
Hence, if every level of ${\Psh{\calC}}$ is induced by a full subcategory of $\calC$ then ${\Psh{\emin{\calC}} \rightarrow \Psh{\calC}}$ is level~{é} of ${\Psh{\calC} \rightarrow \Set}$. 
\end{lemma}
\begin{proof}
Let $B$ be an object in $\calB$ and let ${f : B\rightarrow C}$ be a map in $\calC$.
By hypothesis, ${f = m e}$ for some monic $m$ and some split-epi $e$. Let $s$ be a section of $e$ and consider  the idempotent ${s e : B \rightarrow B}$. By \cite[Lemma~{C5.2.4}]{elephant}, every map in $\calB$ is monic. In particular, the idempotent ${s e}$ is monic, so $e$ must be an iso and,  hence, $f$ is monic. Therefore, $B$ is minimal in $\calC$. 
\end{proof}

For instance, \cite[Proposition~{4.10}]{KellyLawvere89}  shows that for finite $\calC$, every subtopos of ${\Psh{\calC}}$ is induced by a full subcategory of $\calC$. We need something more general. So, instead of finiteness, we will rely on an ascending chain condition for subobjects; but first we prove an auxiliary result, assuming that the reader is familiar with the bijective correspondence \cite[Theorem~{4.4}]{KellyLawvere89} between levels of a presheaf topos and idempotent ideals of the site.

\begin{lemma}\label{LemToBuildAC}
Let $\calC$ be such that every map  factors as a split-epimorphism followed by a monomorphism. 
If $\calI$ is an idempotent ideal in $\calC$ then every ${h \in \calI}$ factors as ${h = i j}$ with ${i, j \in \calI}$ and $i $ monic.
\end{lemma}
\begin{proof}
Because $\calI$ is an idempotent ideal, ${h = f g }$ for some ${f, g \in \calI}$. 
By our factorization hypothesis, ${f = i r}$ with $i$ monic and $r$ split epi.
If we let $s$ be a section of $r$ then ${i = i r s = f s \in \calI}$ because $\calI$ is a right ideal.
Also, ${j = r g \in\calI}$ because $\calI$ is a left ideal. Since ${i j = i r g = f g = h}$,  the result follows.
\end{proof}

We next formulate the chain condition that we need.

\begin{definition}\label{DefACC}
A category $\calC$ satisfies the {\em ascending chain condition (ACC)} if for every  commutative diagram as below
\[ \xymatrix{
C_0 \ar@(d,l)[drr]|-{i_0} \ar[r]^-{m_0} & C_1 \ar[dr]|-{i_1} \ar[r]^-{m_1} & \ldots \ar[r] & \ar[ld]|-{i_{n}} C_n  \ar[r]^-{m_n} & C_{n+1} \ar@(d,r)[dll]|-{i_{n+1}}  \ar[r]^-{m_{n+1}} & \ldots \\
                  &                   &   C    &
}\]
with ${n\in \Nat}$ and where every map is monic, there is a ${t \in \Nat}$ such that ${m_t}$ is an isomorphism. 
\end{definition}

For instance, if every object in $\calC$ has a finite number of subobjects then the ACC  holds.

The next result extends  \cite[Proposition~{3.2}]{MenniPAMS} which shows that: if every map in the small $\calC$ factors as a split epi followed by
a monomorphism then the mono-cartesian idempotent ideals in $\calC$ are in bijective correspondence with the full subcategories of $\calC$ that are closed under subobjects.

\begin{proposition}\label{PropLevelsAreSubcats}
If $\calC$ is such that:
\begin{enumerate}
\item every map factors as a split epi followed by a monic and 
\item satisfies the ACC,
\end{enumerate}
then every level of ${\Psh{\calC}}$ is induced by a full subcategory of $\calC$. 
Hence, ${\Psh{\emin{\calC}} \rightarrow \Psh{\calC}}$ is level~{é} of ${\Psh{\calC} \rightarrow \Set}$.  
\end{proposition}
\begin{proof}
A level of $\Psh{\calC}$ is determined by an idempotent ideal ${\calI}$ of $\calC$.
Define ${\calB \rightarrow \calC}$ to be the full subcategory consisting of those objects whose corresponding identity map is in $\calI$, and let ${\calJ}$ be the idempotent ideal of $\calC$ consisting of the maps that factor through some object in ${\calB\rightarrow \calC}$. Clearly ${\calJ \subseteq \calI}$, so we must prove that ${\calI \subseteq \calJ}$. Let ${f : D \rightarrow C}$ in $\calI$. By Lemma~\ref{LemToBuildAC}, ${f = i_0 j}$ with ${i_0 , j \in \calI}$ and $i_0$ monic. 
We may then apply the same lemma to $i_0$ and iterate in order to build a diagram 
\[ \xymatrix{
D \ar@(d,ld)[rrrd]_-{f} \ar[r]^-{j} & C_0 \ar@(d,l)[drr]|-{i_0} \ar[r]^-{m_0} & C_1 \ar[dr]|-{i_1} \ar[r]^-{m_1} & \ldots \ar[r] & \ar[ld]|-{i_{n}} C_n  \ar[r]^-{m_n} & C_{n+1} \ar@(d,r)[dll]|-{i_{n+1}}  \ar[r]^-{m_{n+1}} & \ldots \\
&                   &                   &   C    &
}\]
where $m_k , i_k \in \calI$ are monic  for every ${k\in \Nat}$.
By the ACC, there is a $k$ such that ${m_k}$ is an iso.
It follows that $C_k$ is in ${\calB\rightarrow \calC}$ and, hence, $f$ is in $\calJ$.

Lemma~\ref{LemPresheafSubetenduesAreBelowTheCanonicalOne} completes the proof.
\end{proof}

We will need the following remark about slicing.

\begin{lemma}\label{LemStabilityOfACC} Let $\calC$ be a small category.
If $\calC$ satisfies the hypotheses of Proposition~\ref{PropLevelsAreSubcats} then so does ${\calC/X}$ for every $X$ in $\Psh{\calC}$.
\end{lemma}
\begin{proof}
Let ${f : (x, C) \rightarrow (y, D)}$ be a map in $\calC/X$.
Let ${f = m e}$ in $\calC$ with ${m : B \rightarrow D}$ monic and $e$ split epi, say, with section ${s : B \rightarrow C}$.
Then $f$ factors as 
\[\xymatrix{
(x, C) \ar[r]^-e & (B, y\cdot m) \ar[r]^-m & (y, D)
}\]
and, since ${x \cdot s = y \cdot m \cdot e \cdot s = y\cdot m}$, the map ${s : (y\cdot m, B) \rightarrow (x, C)}$ is a section of the map ${e :(x, C) \rightarrow (y\cdot m, B)}$. Also, since ${\calC/X \rightarrow \calC}$ is faithful ${m : (y\cdot m, B) \rightarrow  (y, D)}$ is monic. Altogether, every map in $\calC/X$ factors as a split epi followed by a monomorphism.

As $\calC$ preserves monomorphisms and reflects isomorphisms, the ACC holds in $\calC/X$.
\end{proof}

Recall that for any $X$ in $\Psh{\calC}$, the inclusion ${\emin{(\calC/X)} \rightarrow \calC/X}$ induces a  level that we denoted by 
${ \rmD X \rightarrow  \Psh{\calC}/X }$  in Section~\ref{SecD}.

\begin{theorem}\label{ThmACC}
Let $\calC$ be a small category where every map  factors as a split epi followed by a monic.
If $\calC$ satisfies the ACC then, for any object $X$ in $\Psh{\calC}$, the composite ${\Psh{\calC}/X \rightarrow \Psh{\calC} \rightarrow \Set}$ has a level~{é}, and it coincides with ${\rmD X \rightarrow  \Psh{\calC}/X}$.
\end{theorem}
\begin{proof}
By  Lemma~\ref{LemStabilityOfACC},  we may apply Proposition~\ref{PropLevelsAreSubcats} to complete the proof.
\end{proof}

It is well-known that every map in  $\Delta$ factors as a split epi followed by a split monic. 
Since every object has a finite number of subobjects,  $\Delta$ satisfies the ACC.
Similarly, the  category $\mathbb{F}$ of finite non-empty sets  also satisfies the hypotheses of Theorem~\ref{ThmACC}.

It would be instructive to compare the  toposes ${\aE X}$ described here with the ${\petit(X)}$ entertained in \cite{Lawvere91}; also with the calibrations described in \cite{Johnstone2012} and, in particular, with the petit toposes  ${\calS(X)}$ of bundles with discrete fibers studied in \cite{Menni2022}. They will be different in general because, for a pre-cohesive ${p : \calE \rightarrow \calS}$, the hyperconnected geometric morphism ${\calE/X \rightarrow \calS(X)}$ is not always local; but the fact that ${\aE Y}$ and ${\calS(Y)}$ coincide for the reflexive graph $Y$ consisting of a single loop is suggestive; see  Example~\ref{ExNonLocalicNew}.

\section{Height and dimension}
\label{SecBeingAndBecoming}

Fix a small category $\calC$. For any object $X$ in $\Psh{\calC}$, we study in this section the bounded depth formulas that hold in  ${\rmD X}$ and how that relates to the dimension of $X$.

\begin{definition}\label{DefHeight}
For any ${n \in \Nat}$, we say that an object $C$  has {\em height below $n$} if every sequence 
\[ \xymatrix{
 \cdot \ar[r] & \cdot \ar[r] & \cdots  \ar[r]& \cdot \ar[r] & C
}\]
of $n+1$ monomorphisms  contains an isomorphism. In this case we write ${\height C \leq n}$.
\end{definition}

For a  concrete example, we consider the case of the site for simplicial sets.

\begin{lemma}\label{LemHeights} 
For any ${[m] \in \Delta}$ and ${n \in \Nat}$,  ${\height [m] \leq n}$ if and only ${m \leq n}$. 
\end{lemma}
\begin{proof}
The inequality    ${\height [m] \leq n}$ holds if and only if every sequence 
\[ \xymatrix{
 \cdot \ar[r] & \cdot \ar[r] & \cdots \ar[r] & [m]
}\]
of $n+1$ injections in $\Delta$  contains an isomorphism. 
This holds if and only if $[m]$ contains at most ${n + 1}$ elements. Equivalently, ${m + 1 \leq n + 1}$.
\end{proof}

The analogue for $\mathbb{F}$, the category of non-empty finite sets, holds similarly.

\begin{lemma}\label{LemFixedHeightClosedUnderMs}
If ${B \rightarrow C}$ is monic and ${\height C \leq n \in \Nat}$ then ${\height B \leq n}$.
\end{lemma}
\begin{proof}
Consider a chain  
\[ \xymatrix{
 \cdot \ar[r] & \cdot \ar[r] & \cdots \ar[r] & A \ar[r]^-a & B &&  \cdot \ar[r] & \cdot \ar[r] & \cdots \ar[r] & A \ar[r]^-{ba} & C
}\]
 of $n+1$ monomorphisms. If we let $b$ be the monic  ${B \rightarrow C}$ in the statement then we can consider the chain on the right above.
 Since ${\height C \leq n}$ by hypothesis, one of the maps in the second chain is an iso.
That is, either one of the maps `before $A$'  is an iso, or  ${b a : B \rightarrow C}$ is an iso.
In the latter case, as $b$ is monic,  $b$ is also an iso; but then $a$ is also an iso. 
 In either case, one of the maps in the original chain must be an iso so ${\height B \leq n}$.
\end{proof}

The following result is one of the key technical observations of the paper.

\begin{proposition}\label{PropHeightsAndIBDnewNew} 
Let $X$ be an object of  ${\Psh{\calC}}$ and let ${n\in \Nat}$.
If ${\height C \leq n}$ for every   ${(x, C)}$ in  $\emin{(\calC/X)}$  then,  ${\rmD X}$ satisfies ${\IBD{n}}$.
If, moreover,  $X$ is strongly regular then the converse holds.
\end{proposition}
\begin{proof}
By Proposition~\ref{PropCharIBD}, ${\rmD X =  \Psh{\emin{(\calC/X)}}}$ satisfies  ${\IBD{n}}$ if, and only if, for every chain
\[ \xymatrix{
(x_n, C_n) \ar[r] & \cdot \ar[r] & \cdots \ar[r] &  (x_1, C_1) \ar[r] & (x_0, C_0)  \ar[r] & (x, C)
}\]
of ${n + 1}$ maps in $\emin{(\calC/X)}$, one of the maps  is an isomorphism. 
Notice that since the objects are minimal, all the maps in such a chain are monic in $\calC/X$, and therefore monic in $\calC$.
So, if  ${\height C \leq n}$,  then one of the monic maps in the chain in $\calC$ is an isomorphism.
As ${\calC/X \rightarrow \calC}$ reflects isomorphisms, the corresponding map in the chain in $\calC/X$ is an isomorphism.

For the converse assume that $X$ is strongly regular and    that  ${\rmD X}$ satisfies ${\IBD{n}}$.
Let ${(x, C)}$ be minimal in $\calC/X$ and let 
\[ \xymatrix{
C_n  \ar[r]  & \cdots  \ar[r]& C_0 \ar[r] & C
}\]
be a sequence of $n+1$ monomorphisms in $\calC$. 
Such a sequence determines in the obvious way a sequence
\[ \xymatrix{
(x_n, C_n) \ar[r]  & \cdots \ar[r]  & (x_0, C_0)  \ar[r] & (x, C)
}\]
 in $\calC/X$. 
 Since $X$ is strongly regular by hypothesis, this sequence is in ${\emin{(\calC/X)}}$  and, therefore, Proposition~\ref{PropCharIBD} implies that 
one of the maps in the sequence is an isomorphism, showing that ${\height C \leq n}$.
\end{proof}

Any principal Lawvere-Tierney topology (and so, in particular, any essential subtopos) determines an associated idempotent {\em principal comonad} \cite{MenniPAMS}. As an instance of this, let ${\calC_n \rightarrow \calC}$ be the full subcategory determined by the objects that have height below $n$. The inclusion determines a level that we denote by ${n : \Psh{\calC_n} \rightarrow \Psh{\calC}}$.  
The counit of ${n_! \dashv n^*}$ may be factored as an epi followed by a monic, and the monic part is the counit of the principal comonad determined by level $n$. This principal comonad will be denoted by  ${\Sk_n}$.  
For any object $X$ in $\Psh{\calC}$, the (monic) counit  ${\Sk_n X \rightarrow X}$ will be called the {\em $n$-skeleton}  of $X$. By \cite[Theorem~{2.6}]{MenniPAMS}, for each object  $D$ in $\calC$, ${(\Sk_n X)D \rightarrow X D}$ consists of  those ${x \in X D}$ such that there is a map ${f : D \rightarrow E}$ in $\calC$ that factors through $\calC_n$ and a ${y \in  X E}$ such that ${x = y \cdot f}$.

Assuming that $\calC$ is equipped with a factorization system it is possible to strengthen the description of the principal comonads above.
We do not work in maximum generality. Recall that a map is a {\em strong epimorphism} if it is left orthogonal to every monomorphism.

\begin{lemma}\label{LemCharSkeletal} 
Assume that every map in $\calC$  factors as a strong epi followed by a monic.
For any  $X$ in $\Psh{\calC}$,  $D$ in $\calC$,  and ${x \in X D}$, ${x \in (\Sk_n X)D \rightarrow X D}$ if, and only if,  there exists an strong epi  ${g : D \rightarrow C}$ such that ${\height C \leq n}$ and a ${z \in X C}$ such that ${z \cdot g = x}$.
\end{lemma}
\begin{proof}
One direction is immediate. 
For the other assume that ${x \in (\Sk_n X)D \rightarrow X D}$ so that there is a map ${f : D \rightarrow E}$ in $\calC$ that factors through $\calC_n$ and a ${y \in  X E}$ such that ${x = y \cdot f}$. Let ${f = m h}$ with ${h : D \rightarrow C}$ such $C$ in $\calC_n$ and, using the factorization hypothesis,  let ${h = b g}$ with ${b : B \rightarrow C}$ monic  and $g$ strongly epi.
Then 
\[ x = y \cdot f = y \cdot (m h) = (y \cdot m) \cdot (bg) = (y\cdot m  \cdot b) \cdot g \]
with ${g : D \rightarrow B}$ strongly epi, and ${\height B \leq n}$ by Lemma~\ref{LemFixedHeightClosedUnderMs}.
\end{proof}

We say  that $X$ is {\em $n$-skeletal} if its $n$-skeleton is an iso. In this case we write ${\dim X \leq n}$.

\begin{proposition}\label{PropDimImpliesIBD} 
Let $\calC$ be a small category where every map factors as a strong epi  followed by a monic. 
For any ${n \in \Nat}$ and any object $X$ in ${\Psh{\calC}}$, if ${\dim X \leq n}$ then ${\IBD{n}}$ holds in ${\rmD X}$.
\end{proposition}
\begin{proof}
Let ${(x, D)}$ be an object in ${\emin{(\calC/X)}}$.
Since ${\dim X \leq n}$ by hypothesis, Lemma~\ref{LemCharSkeletal} implies the existence of  a strong epimorphism ${g : D \rightarrow C}$ such that ${\height C \leq n}$ and a ${z \in X C}$ such that ${z \cdot g = x}$.
Since ${(x, D)}$ is minimal, $g$ must be an isomorphism. So ${\height D = \height C \leq n}$ by Lemma~\ref{LemFixedHeightClosedUnderMs} and hence, Proposition~\ref{PropHeightsAndIBDnewNew}  completes the proof.
\end{proof}

\section{Well-foundedness}
\label{SecWF}

In this section we introduce additional hypothesis in order to prove a converse of Proposition~\ref{PropDimImpliesIBD} which becomes the main result of the paper. Fix a small category $\calC$.

\begin{definition}\label{DefWellFounded}
We say that $\calC$ is {\em well-founded} if every chain
\[ \xymatrix{
C_0  \ar[r]^-{e_0} & C_1 \ar[r]^{e_1} &  \cdots \ar[r] & C_n \ar[r]^-{e_{n}} & \cdots 
}\]
of strong epimorphisms is eventually constant (in the sense that there is a ${k \in\Nat}$ such that $e_n$ is an
iso for all ${n \geq k}$ in ${\Nat}$).
\end{definition}

The following result is a   weak form, in a more general context, of the `Eilenberg-Zilber Lemma' for simplicial sets \cite[Proposition~{II.3.1}]{GabrielZisman}. See also   \cite[Lemma~{3.9}]{FioreMenni2005}.

\begin{lemma}\label{LemMinimallyEngendered}
If $\calC$ is  well-founded and every map in $\calC$ factors as a strong epi followed by a monic then, for every $X$ in $\Psh{\calC}$, every $C$ in $\calC$ and every ${x \in X C}$ there exists a strong epi ${e : C \rightarrow D}$ in $\calC$ and a minimal object ${(y, D)}$ in ${\calC/X}$ such that ${x = y \cdot e}$.    
\end{lemma}
\begin{proof}
Assume that there are no  map $e$ and a minimal object ${(y, D)}$ as in the statement.
Then  there is a map ${f_0 : (x, C) \rightarrow (y_1, D_1)}$ in ${\calC/X}$ with non-monic $f$ in $\calC$.
Let ${f_0 = m_0 e_0}$ with strongly epi, and necessarily not iso, ${e_0 : C \rightarrow C_1}$ and monic $m_0$.
Let ${x_1 = y_1 \cdot m_0}$ and consider the resulting map ${e_0 :(x, C) \rightarrow (x_1 , C_1)}$.
Since ${(x_1, C_1)}$ is not minimal we can construct in the same way an ${e_1 : (x_1, C_1) \rightarrow (x_2 , C_2)}$ with strongly epi (and non-iso)  $e_1$ in $\calC$.
Iterating we produce a sequence
\[ \xymatrix{
C_0  \ar[r]^-{e_0} & C_1 \ar[r]^{e_1} &  \cdots \ar[r] & C_n \ar[r]^-{e_{n}} & \cdots 
}\]
with $e_i$ strongly epi and non-iso for every ${i \in \Nat}$, contradicting well-foundedness.
\end{proof}

We may now prove a partial converse to Proposition~\ref{PropDimImpliesIBD}.

\begin{proposition}\label{PropIBDImpliesDim}
Let $\calC$ be  well-founded   and such that every map factors as a strong epi followed by a monic.
For any ${n \in \Nat}$ and  strongly regular $X$ in ${\Psh{\calC}}$,  if  ${\IBD{n}}$ holds in ${\rmD X}$ then ${\dim X \leq n}$. 
\end{proposition}
\begin{proof}
Let ${(x, C)}$ be an object of ${\calC/X}$.
By Lemma~\ref{LemMinimallyEngendered} there is a map ${e : (x, C) \rightarrow (y, D)}$ with $e$ strongly epi and  ${(y, D)}$ minimal in ${\calC/X}$. 
Proposition~\ref{PropHeightsAndIBDnewNew} implies that  ${\height D \leq n}$.
Since this holds for every ${(x, C)}$, Lemma~\ref{LemCharSkeletal} implies that ${\dim X \leq n}$.
\end{proof}

Combining Propositions~\ref{PropDimImpliesIBD}  and~\ref{PropIBDImpliesDim} we obtain the main result of the paper.

\begin{theorem}\label{ThmDimIffIBD}
Let $\calC$ be small and well-founded category where  every map factors as a strong epimorphism followed by a monomorphism.
Let $X$ in ${\Psh{\calC}}$ be strongly regular.
For any ${n \in \Nat}$, ${\dim X \leq n}$ if and only if ${\IBD{n}}$ holds in ${\rmD X}$. 
\end{theorem}

The category ${\Delta}$  satisfies the hypothesis of Theorems~\ref{ThmACC} and~\ref{ThmDimIffIBD}.
Also, it follows from Lemma~\ref{LemHeights} that, for any ${X \in \Psh{\Delta}}$ and ${n \in \Nat}$, ${\dim X\leq n}$ in the sense of Theorem~\ref{ThmDimIffIBD} if and only if ${\dim X \leq n}$ in the usual sense \cite{GabrielZisman}. In other words, $n$-skeletal in the sense of Section~\ref{SecBeingAndBecoming} is equivalent to its usual meaning in the context of simplicial sets.

\begin{corollary}
For any strongly regular $X$ in $\Psh{\Delta}$ and ${n\in\Nat}$, ${\dim X \leq n}$ if, and only if, $\IBD{n}$ holds in ${\aE X}$.
\end{corollary}

The analogous result holds for $\mathbb{F}$ instead of $\Delta$.

We would like to have a rigorous comparison with  \cite{MarraEtAl2018}.
One way to proceed would be as follows. 
First, construct a topos embedding the category of polyhedra used in \cite{MarraEtAl2018}.
For instance, one could use the (non presheaf) PL-topos described by Marra in his  invited talk at the conference {\em Topology, Algebra and Categories in Logic 2019}.
Second, derive an informative  notion of petit topos in that context.
It is to be expected that for a polyhedron $P$, seen as an object in the topos, the associated petit topos could coincide with  the topos of sheaves on the topological space of open subpolyhedra of $P$ used in \cite{MarraEtAl2018}.
The unit interval in the PL-topos is a model of the theory classified by $\Psh{\Delta}$ so it determines a geometric morphism from the former topos to the latter which could be used, not only to compare their homotopy theories as in \cite{MarmolejoMenni2017}, but maybe also to compare the respective notions of petit topos.

A more direct way to proceed is as follows. For a polyhedron $P$, let ${P'}$ be the simplicial set that sends ${[n]}$ to the set of affine maps ${\Delta^n \rightarrow P}$, where $\Delta^n$ is the standard $n$-simplex. A minimal element of ${P'}$ is an affine map ${f : \Delta^n \rightarrow P}$ such that for any factorization ${f = g h' : \Delta^n \rightarrow \Delta^m \rightarrow P}$ with ${g : \Delta^m \rightarrow P}$ affine and ${h' : \Delta^n \rightarrow \Delta^m}$ induced by a monotone ${h : [n] \rightarrow [m]}$, $h$ is injective. So the existence of a minimal ${f : \Delta^n \rightarrow P}$ implies that ${\dim P \geq n}$.
In other words, it seems that ${\rmD(P')}$ has enough information to determine the dimension of $P$ but we will not pursue the idea here.

To a certain extent,  some of the results above  also apply to classical algebraic geometry, via the not-quite-combinatorial  Gaeta topos determined by the field $\mathbb{C}$ which is well-known to be a pre-cohesive presheaf topos over $\Set$. 
Among other reasons, this example would be interesting because, in contrast with the cases of ${\Psh{\Delta}}$ and ${\Psh{\mathbb{F}}}$ emphasized in the present paper, the Gaeta topos for $\mathbb{C}$  has non-trivial level~$\epsilon$. See \cite{MMlevelEpsilon} and also \cite[Example~{6.7}]{MenniPAMS}. 

Also the more recent Rig Geometry is a possible realm of application. Indeed, as shown in \cite{Menni2021a}, the Gaeta toposes for the theories of 2-rigs and integral rigs are  presheaf toposes; and the same holds for the theory of positive rigs  \cite{Menni+2024+Postive+Rigs}.
The fact that the sites for these toposes have epi/regular-mono factorizations suggest that it may be interesting to consider a notion of minimal object relative to a factorization system as in \cite{FioreMenni2005}.

\end{document}